\documentclass[12pt]{article}
\begin{document}
\title{A Sieve for Cousin Primes}
\author{H. J. Weber\\Department of Physics\\
University of Virginia\\Charlottesville, VA 22904\\USA}
\maketitle
\begin{abstract}
A sieve is constructed for twin primes at distance $4,$ which 
are of the form $3(2m+1)\pm 2,$ and are characterized by 
their twin-4 rank $2m+1.$ It does not suffer from the parity 
problem. Non-rank numbers are identified and counted using odd 
primes $p\geq 5.$ Twin-4 ranks and non-ranks make up the set 
of odd numbers. Regularities of non-ranks allow gathering 
information on them to obtain a Legendre-type sum for the 
number of twin-4 ranks. Due to considerable cancellations 
in it, the asymptotic law of its main term has the expected 
form and magnitude of its coefficient.            
\end{abstract}
\vspace{3ex}
\leftline{MSC: 11A41, 11N05}
\leftline{Keywords: Twin-4 rank, non-ranks, sieve} 


\section{Introduction}

Our knowledge of twin primes comes mostly from sieve 
methods~\cite{hr},\cite{rm},\cite{hri},\cite{fi}. 

In Ref.~\cite{adhjw} a sieve is developed specifically 
for ordinary twin primes. These methods are applied here 
to cousin primes, where arithmetical details are rather  
different. However, there are also similarities because 
both their distances, $2$ and $4,$ have no odd prime 
divisor.  

Prime numbers $p\geq 5$ are well known to be of the form~\cite{hw} 
$6m\pm 1.$ An ordinary twin prime occurs when both $6m\pm 1$ are 
prime. Twin primes at distance $4$ can be written similarly as 
$6m+1, 6(m+1)-1$ or $3(2m+1)\pm 2,$ being in class II of a  
classification~\cite{hjw} of all twin primes, whereas ordinary 
twins lie in class I being of the form $2(3m)\pm 1$.  

{\bf Definition~1.1.} If $3(2m+1)\pm 2$ is a twin prime pair for 
some odd $2m+1$, then $2m+1$ is called its {\it twin-4 rank}. An 
odd number $2m+1$ is a {\it non-rank} if $3(2m+1)\pm 2$ are not 
both prime. Odd positive integers $\geq 3$ consist of twin-4- 
and non-ranks only. Even numbers are not considered in the 
following because $3(2m)\pm 2$ are never primes. Also, since 
$2, 3$ are not of the form $6m\pm 1,$ they are excluded as 
primes in the following. Also, we ignore the {\it special} 
cousin prime~\cite{hjw} $5\pm 2=(3,7).$  

{\bf Example~1.2.} Twin-4 ranks are $3, 5, 7, 13, 15,\ldots .$ 
Non-ranks are $9, 11, 17,\\19,\ldots .$  

The odd numbers $\geq 3$ form the base set of this pair sieve; 
it is partitioned into twin-4 and non-rank sets. Only non-ranks 
have sufficient regularity and abundance allowing us to gather 
enough information on them to draw inferences on the minimal 
number of twin-4 ranks needed to account for all odd numbers 
$\geq 3$. Therefore, our main focus is on non-ranks, their 
symmetries and abundance.      

In Sect.~2 the twin-4 prime sieve is constructed based on non-ranks. 
In Sect.~3 non-ranks are identified in terms of their main properties 
and then, in Sect.~4, they are counted. In Sect.~5 twin-4 ranks are 
isolated and counted. Conclusions are summarized and discussed in 
Sect.~6.   

\section{Twin Ranks, Non-Ranks and Sieve}

It is our goal here to construct a twin-4 prime sieve. To this end, 
we need the following arithmetical function. 

{\bf Definition~2.1.} Let $x$ be real. Then $N(x)$ is the integer 
nearest to $x.$ The ambiguity for $x=n+\frac{1}{2}$ with integral 
$n$ will not arise in the following.    

{\bf Lemma~2.2.} {\it Let $p\geq 5$ be prime. Then} 
\begin{eqnarray}
N(\frac{p}{6})=\{\begin{array}{ll}\frac{p-1}{6},~\rm{if}~p\equiv 
1\pmod{6};\\\frac{p+1}{6},~\rm{if}~p\equiv -1\pmod{6}.\\
\end{array}
\end{eqnarray}
{\it Proof.} This is obvious from Def.~2.1 by substituting 
$p=6m\pm 1.~\diamond$ 

{\bf Corollary~2.3.} {\it If $p\equiv -1\pmod{6}$ is prime and 
$p-4$ is prime, then $\frac{p+1}{3}-1$ is a twin-4 rank. If 
$p\equiv 1\pmod{6}$ and $p+4$ is prime, then $\frac{p-1}{3}+1$ 
is a twin-4 rank.}  

{\it Proof.} This follows from Def.~1.1, as $3(\frac{p+1}{3}
-1)\pm 2=(p, p-4),$ and $3(\frac{p-1}{3}+1)\pm 2=(p, p+4)$ 
in the last case.$~\diamond$  

{\bf Example~2.4.} This is the case for $p=7,13,19,\ldots$ as 
well as for $p=11,17,23,\dots$ but not for $p=5, 29, 31,\dots .$  

{\bf Lemma~2.5} {\it Let $p\geq 5$ be prime. Then all odd  
numbers 
\begin{eqnarray}\nonumber 
&&2k(n,p)^++1=(2n+1)p+4N(\frac{p}{6}),~n=0,1,2,\ldots\\
&&2k(n,p)^-+1=(2n+1)p-4N(\frac{p}{6}),~n=1,2,\ldots 
\end{eqnarray}
are non-ranks. There are $2=2^{\nu(p)}$ (single) non-rank 
progressions to the prime $p,$ where $\nu(n)$ counts the 
number of different prime divisors of $n.$  
  
(a) If $p\equiv 1\pmod{6}$ the non-rank $2k(n,p)^++1$ 
generates the pair}  
\begin{eqnarray}
3(2k(n,p)^++1)\pm 2=([3(2n+1)+2]p-4,[3(2n+1)+2]p),
\label{p++}
\end{eqnarray} 
{\it and the non-rank $2k(n,p)^-+1$ the pair}  
\begin{eqnarray}
3(2k(n,p)^-+1)\pm 2=([3(2n+1)-2]p,[3(2n+1)-2]p+4).
\label{p+-}
\end{eqnarray} 
 
{\it (b) If $p\equiv -1\pmod{6}$ the non-rank 
$2k(n,p)^++1$ generates the pair}  
\begin{eqnarray}
3(2k(n,p)^++1)\pm 2=([3(2n+1)+2]p-4,[3(2n+1)+2]p),
\label{p-+}
\end{eqnarray} 
{\it and the non-rank $2k(n,p)^-+1$ the pair} 
\begin{eqnarray}
3(2k(n,p)^-+1)\pm 2=([3(2n+1)-2]p-4,[3(2n+1)-2]p).
\label{p--}
\end{eqnarray}
{\it All pairs contain a composite number.} 

For $n=0$ and $p\equiv \pm 1\pmod{6},~2k^-+1
\to\frac{p\pm 2}{3}$ are the twin-4 ranks of 
Cor.~2.3.

Clearly, all these non-ranks are symmetrically distributed 
at equal distances $4N(p/6)$ from odd multiples of each 
prime $p\geq 5.$ Twin-4 ranks and some non-ranks for $n=0$ 
are the subject of Cor.~2.3 and Example~1.2. 
  
{\it Proof.} Let $p\equiv 1\pmod{6}$ be prime and $n\geq 0$ an 
integer. Then $2k(n,p)^++1=(2n+1)p+4\frac{p-1}{6}$ by Lemma~2.2 
and $3(2k^++1)$ is sandwiched by the pair in Eq.~(\ref{p++}) 
which contains a composite number. Hence $2k(n,p)^++1$ is a 
non-rank. For $n>0,$ the same happens in Eq.~(\ref{p+-}), so 
$2k^-+1$ is a non-rank. 

If $p\equiv -1\pmod{6}$ and prime, then $2k(n,p)^++1=(2n+1)p+ 
4\frac{p+1}{6}$ by Lemma~2.2 and $3(2k^++1)$ leads to the pair 
in Eq.~(\ref{p-+}) which contains a composite number again. 
For $n>0,$ the same happens in Eq.~(\ref{p--}), so $2k^-+1$ is 
a non-rank.~$\diamond$

The $2k(n,p)^{\pm}+1$ yield pairs $3(2k^{\pm}+1)\pm 2$ with one 
or two composite entries that are twin-4 prime analogs of multiples 
$np,~n>1$ of a prime $p$ in Eratosthenes' prime sieve~\cite{hw}. 

The converse of Lemma~2.5 holds, i.e. non-ranks are organized 
in arithmetic progressions by prime numbers $\geq 5.$    
 
{\bf Lemma~2.6.} {\it If $2k+1>3$ is a non-rank, there is a prime 
$p\geq 5$ and a non-negative odd integer $2\kappa\pm 1$ so that 
$2k+1=2k(\kappa,p)^++1$ or} $2k(\kappa,p)^-+1.$ 

{\it Proof.} Let $6k+5=3(2k+1)+2$ be composite. Then 
$3(2k+1)+2=6k+5\neq 2^\mu 3^\nu,~\mu+\nu\geq 1$ obviously. 
Let $6k+5=p\cdot K\equiv -1\pmod{6},$ where $p\geq 5$ is the 
smallest prime divisor. If $p=6m+1,$ then $K=6\kappa-1$ and  
\begin{eqnarray}
6k+5=6^2 m\kappa+6(\kappa-m)-1,~k+1=6m\kappa+\kappa-m=p\kappa
-\frac{p-1}{6}
\end{eqnarray}
and 
\begin{eqnarray}
2k+1=2p\kappa-\frac{p-1}{3}-1=(2\kappa-1)p+4\frac{p-1}{6},
\end{eqnarray}       
q.e.d. If $p=6m-1,$ then $K=6\kappa+1$ and  
\begin{eqnarray}
6k+5=6^2 m\kappa+6(m-\kappa)-1,~k+1=6m\kappa+m-\kappa=p\kappa
+\frac{p+1}{6}
\end{eqnarray}
and so 
\begin{eqnarray}
2k+1=2p\kappa+\frac{p+1}{3}-1=(2\kappa+1)p-4\frac{p+1}{6},
\end{eqnarray}       
q.e.d. Now let $6k+1=pK$ with $p\equiv 1\pmod{6}.$ Then 
$K=6\kappa+1$ and, therefore, 
\begin{eqnarray}
k=6m\kappa+m+\kappa=p\kappa+\frac{p-1}{6},~
2k+1=(2\kappa+1)p-4\frac{p-1}{6}, 
\end{eqnarray}
q.e.d. Finally, if $p\equiv -1\pmod{6}$ then $K=6\kappa-1$ 
and 
\begin{eqnarray}\nonumber
&&6k+1=(6m-1)(6\kappa-1)=6^2m\kappa-6(m+\kappa)+1,\\
&&k=6m\kappa-(m+\kappa)=p\kappa-\frac{p+1}{6}.
\end{eqnarray}
Hence
\begin{eqnarray}
2k+1=(2\kappa-1)p+4\frac{p+1}{6},
\end{eqnarray}
q.e.d.$~\diamond$ 
 
Even multiples of prime numbers $p\geq 5$ in Lemma~2.5, 
e.g. $2np+1\pm 2\alpha N(p/6)$ for appropriate $\alpha,$ 
are accounted for in Lemma~2.6 as non-ranks to some prime 
$p'\geq 5$, which demonstrates the cornerstone role 
Lemma~2.6 plays for the sieve. 

{\bf Theorem~2.7. (Cousin Prime Sieve)} {\it Let} ${\cal 
P}=\{(2n+1\geq 3,2n+5): n\geq 0,~\rm{integral}\}$ {\it be 
the set of all pairs with entries $\geq 3$ of natural 
numbers at distance $4.$ Upon striking all pairs identified 
by non-ranks of Lemma~2.5, only (and all) twin-4 prime pairs 
are left.} 

Clearly, this sieve is not subject to the parity problem.    
 
{\it Proof.} Obviously, we need to consider only the subset 
${\cal P}_0=\{(6m+1,\\6m+5): m\geq 0,~\rm{integral}\}
\subset{\cal P}.$ For $2m+1\geq 3$ divide $3(2m+1)\pm 2$ 
by all primes $p<\sqrt{6m+1}.$ Then $2m+1$ is a non-rank if 
there is a prime $p$ such that $(6m+5)/p$ or $(6m+1)/p$ 
(or both) is integral. For all such $m,~2m+1$ is struck 
from the set of odd positive integers. Then all remaining 
odd integers are twin-4 ranks.~$\diamond$    

More concrete steps to construct it will be taken in the 
next section. 
 
\section{Identifying Non-Ranks} 

Here it is our goal to systematically characterize 
and identify non-ranks among odd numbers. 

{\bf Definition~3.1} Let $p\geq 5$ be the minimal 
prime of a non-rank. Then $p$ is its parent prime.  

{\bf Example~3.2.} The non-ranks to parent prime 
$5$ are, by Lemma~2.5,   
\begin{eqnarray}
2k^++1=9,19,29,\dots ;~2k^-+1=11,21,\ldots
\end{eqnarray} 
These $2k^{\pm}+1$ form the set ${\cal A}_5^-=
\{5(2n+1)\pm 4\geq 9:~n\geq 0\}={\cal A}_5.$ Note that 
$5$ is the most effective non-rank generating prime 
number. If it were excluded like $3$ then many numbers, 
such as $9, 19,\ldots,$ would be missed as non-ranks. 

{\bf Proposition~3.3.} {\it The arithmetic progressions} 
$3\cdot 5(2n+1)\pm 2,3[5(2n+1)+2]\pm 2,3[5(2n+1)-2]\pm 2$ 
{\it contain all twin-4 prime pairs.} 

Prop.~3.3 is the first step of the twin-4-prime sieve.
Let ${\cal C}_5=\{0,2,8\}$ be the set of non-negative 
constants $c$ in $5(2n+1)+c$ in Prop.~3.3. 

{\it Proof.} From $\{3\cdot 5(2n+1)\pm 2,~3[5(2n+1)+2]
\pm 2, 3[5(2n+1)-2]\pm 2, 3[5(2n+1)+4]\pm 2, 3[5(2n+1)
-4\pm 2,~n\geq 0\}$ we strike all pairs $\{3\cdot 
5(2n+1)+12\pm 2,~3\cdot 5(2n+1)-12\pm 2,~n>0\}$ 
resulting from non-ranks of ${\cal A}_5^-.~\diamond$ 

For $p=7,$ we now subtract from the set ${\cal A}_7^+
=\{7(2n+1)\pm 4\geq 9:~n\geq 0\}$ of non-ranks to 
$7$ those to $p=5.$ The remaining set ${\cal A}_7$ 
comprises the non-ranks to parent prime $p=7.$ 
 
{\bf Lemma~3.4.} {\it The set ${\cal A}_7$ of 
non-ranks to parent prime $p=7$ comprises the 
arithmetic progressions} $\{7\cdot 5(2n+1)+10,
7\cdot 5(2n+1)+18,7\cdot 5(2n+1)+32,7\cdot 
5(2n+1)+38,7\cdot 5(2n+1)+52,7\cdot 5(2n+1)+60\}.$

{\it Proof.} We subtract the common arithmetic 
progressions of ${\cal A}_5^-: \{5(2n+1)\pm 4: 
2n+1\to 7(2n+1), 7(2n+1)\pm 2,7(2n+1)\pm 4,
7(2n+1)\pm 6\}$ from ${\cal A}_7^+: \{7(2n+1)\pm 
4: 2n+1\to 5(2n+1),5(2n+1)\pm 2,5(2n+1)\pm 4\}$ 
to find those listed in Lemma~3.4.  
 
The common arithmetic progressions are 
\begin{eqnarray}\nonumber
&&5\cdot 7(2n+1)\pm 4,\\\nonumber
&&5[7(2n+1)+4]+4=7[5(2n+1)+4]-4,\\
&&5[7(2n+1)-4]-4=7[5(2n+1)-4]+4.~\diamond  
\end{eqnarray} 

Note that these $2^{\nu(5\cdot 7)}=2^2$ 
arithmetic progressions contain all common 
(double) non-ranks of the primes $5,~7.$ 

{\bf Proposition~3.5.} {\it The arithmetic 
progressions} $3[5\cdot 7(2n+1)+c]\pm 2, 
~n\geq 0$ {\it contain all twin-4 pairs $\geq 
103$, where $c\in {\cal C}_7=\{0,2,8,12,20,22,
28,30,40,42,\\48,50,58,62,68\}$.} 

Note that ${\cal C}_5\subset {\cal C}_7,$ but 
this pattern does not continue.   

{\it Proof.} Using Lemma~3.4, we strike from the 
arithmetic progressions of Prop.~3.3 (replacing 
$2n+1\to 7(2n+1),7(2n+1)\pm 2,7(2n+1)\pm 4, 
7(2n+1)\pm 6$) all pairs resulting from 
non-ranks in ${\cal A}_7,$ which are $\{5\cdot 
7(2n+1)+a; a=10,18,32,38,52,60.\}$ This leaves 
the progressions listed above.$~\diamond$  

This is the second step of the sieve. 

In contrast to ordinary twin primes the 
arithmetic function values $N(p'/6),\\N(p/6)$ 
do not suffice to characterize twin-4 primes 
$p'=p+4.$ 

{\bf Theorem~3.6.} {\it Let $p', p$ be primes. 
If $p'\equiv -1\pmod{6}, p\equiv 1\pmod{6}$ and 
$N(\frac{p'}{6})=1+N(\frac{p}{6})$ then} $p'=p+4.$

{\it Proof.} If $p'\equiv -1\pmod{6}, p\equiv 
1\pmod{6}$ then 
\begin{eqnarray}
\frac{p'+1}{6}=1+\frac{p-1}{6}
\end{eqnarray}
is equivalent to $p'=p+4.~\diamond$

{\bf Corollary~3.7.} {\it Let $p'>p\geq 5$ be 
primes such that $N(\frac{p'}{6})=1+N(\frac{p}
{6}).$ Then $p'=p+4$ if $p\equiv 1\pmod{6}$ 
and $p'\equiv -1\pmod{6};$ if $p'\equiv 
1\pmod{6}$ instead then $p'=p+6.$ If 
$p\equiv -1\pmod{6}$ and $p'\equiv 
-1\pmod{6}$ then $p'=p+6;$ if $p'\equiv 
1\pmod{6}$ instead then} $p'=p+8.$
 
{\it Proof.} Let $p\equiv -1\pmod{6}.$ Then 
$1+N(\frac{p}{6})=1+\frac{p+1}{6}.$ If 
$p'\equiv -1\pmod{6}$ then $p'=p+6.$ If  
$p'\equiv 1\pmod{6}$ then $p'=p+8.$ If 
$p\equiv 1\pmod{6}$ then $1+\frac{p-1}{6}=
N(p'/6)$ implies $p'=p+6$ if $p'\equiv 
1\pmod{6};$ if $p'\equiv -1\pmod{6}$ then  
$p'=p+4.~\diamond$ 

{\bf Theorem~3.8.} {\it Let $p\geq 5,~
p'=p+2$ be ordinary prime twins. Then 
$pp'(2n+1)\pm 4\frac{p+1}{6}>0$ for 
$n=0,1,2,\ldots$ and  
\begin{eqnarray}\nonumber
&&p[p'(2n+1)+4\frac{p+1}{6}]+4\frac{p+1}{6}=
p'[p(2n+1)+4\frac{p+1}{6}]-4\frac{p'-1}{6}>0,\\
\nonumber
&&p[p'(2n+1)-4\frac{p+1}{6}]-4\frac{p+1}{6}=
p'[p(2n+1)-4\frac{p+1}{6}]+4\frac{p'-1}{6}>0,\\
&&n=0,1,2,\ldots 
\label{nonr}
\end{eqnarray} 
are their common non-ranks.}   

Note that again there are $4$ arithmetic  
progressions of common or double non-ranks. 

{\it Proof.} Using $N(p'/6)=N(p/6)=\frac{p+1}
{6},$ Eq.~(\ref{nonr}) is readily verified;  
its lhs $\in {\cal A}_p^-$ and rhs $\in 
{\cal A}_{p'}^+$ and $p(p+2)(2n+1)\pm 4\frac{p+1}
{6}\in {\cal A}_p^-,{\cal A}^+_{p'}$.~$\diamond$ 

We now consider systematically common  
non-ranks of pairs of primes generalizing 
Theor.~3.8 to arbitrary prime pairs.  

{\bf Theorem~3.9.} {\it Let $p'>p\geq 5$ 
be primes. (i) If $p'\equiv p\equiv -1\pmod{6},$ 
then $p'=p+6l,~l\geq 1,~ N(\frac{p'}{6})=
N(\frac{p}{6})+l$ and common non-ranks of 
$p',p$ are, for} $n=0, 1,\ldots,$
\begin{eqnarray}
p[p'(2n+1)+2r']\pm 4N(\frac{p}{6})=p'[p(2n+1)+2r]\pm 
4N(\frac{p'}{6})
\label{snr1}
\end{eqnarray}
{\it provided the integers $r, r'$ solve} 
\begin{eqnarray}
(r'-r)p=2l(3r\pm 1),~1-p\leq 2r\leq p-1,~1-p'\leq 
2r'\leq p'-1.
\label{snr11}
\end{eqnarray}
{\it Eq.~(\ref{snr11}) with $3r\pm 1\equiv 0\pmod{p}$ 
on the rhs has a unique solution $r$ that then 
determines} $r'.$ 

{\it If $r, r'$ solve} 
\begin{eqnarray}
(r'-r)p=2l(3r\mp 1)l\mp 2N(\frac{p+1}{3}) 
\label{snr12}
\end{eqnarray}
{\it then the common non-ranks are}
\begin{eqnarray}
p[p'(2n+1)+2r']\pm 4N(\frac{p}{6})=p'[p(2n+1)
+2r]\mp 4N(\frac{p'}{6}).
\label{snr13}
\end{eqnarray}
{\it (ii) If $p'\equiv p\equiv 1\pmod{6},$ 
then $p'=p+6l,~l\geq 1,~N(\frac{p'}{6})=
N(\frac{p}{6})+l,$ and common non-ranks of 
$p',p$ are}
\begin{eqnarray}
p[p'(2n+1)+2r']\pm 4N(\frac{p}{6})=p'[p(2n+1)
+2r]\pm 4N(\frac{p'}{6})
\label{snr2}
\end{eqnarray}
{\it provided $r, r'$ solve}
\begin{eqnarray}
(r'-r)p=2l(3r\pm 1). 
\label{snr22}
\end{eqnarray}
{\it If $r, r'$ solve} 
\begin{eqnarray}
(r'-r)p=2l(3r\mp 1)l\mp 2N(\frac{p-1}{6}) 
\label{snr23}
\end{eqnarray}
{\it then the common non-ranks are}
\begin{eqnarray}
p[p'(2n+1)+2r']\pm 4N(\frac{p}{6})=p'[p(2n+1)
+2r]\mp 4N(\frac{p'}{6}).
\label{snr24}
\end{eqnarray}
{\it (iii) If $p'\equiv 1\pmod{6},~p\equiv 
-1\pmod{6}$ then $p'=p+6l+2,~l\geq 0,
~N(\frac{p'}{6})=N(\frac{p}{6})+l,$ and 
common non-ranks of $p',p$ are}
\begin{eqnarray}
p[p'(2n+1)+2r']\pm 4N(\frac{p}{6})=p'[p(2n+1)
+2r]\pm 4N(\frac{p'}{6})
\label{snr3}
\end{eqnarray}
{\it provided}
\begin{eqnarray}
(r'-r)p=2r(3l+1)\pm 2l. 
\label{snr33}
\end{eqnarray}
{\it If $l=0$ then $r'=r=0$ and Eq.~(\ref{nonr}) 
are solutions (Cor.~3.7).} 

{\it If $r, r'$ solve} 
\begin{eqnarray}
(r'-r)p=2r(3l+1)l\mp 2\left(l+\frac{p+1}{3}
\right),~l\geq 1,  
\label{snr31}
\end{eqnarray}
{\it then the common non-ranks are}
\begin{eqnarray}
p[p'(2n+1)+2r']\pm 4N(\frac{p}{6})=p'[p(2n+1)
+2r]\mp 4N(\frac{p'}{6}).
\label{snr32}
\end{eqnarray}
{\it (iv) If $p'\equiv -1\pmod{6},~p\equiv 
1\pmod{6}$ then $p'=p+6l-2,~l\geq 1,
~N(\frac{p'}{6})=N(\frac{p}{6})+l,$ and
common non-ranks of $p',p$ are}
\begin{eqnarray}
p[p'(2n+1)+2r']\pm 4N(\frac{p}{6})=p'[p(2n+1)
+2r]\pm 4N(\frac{p'}{6})
\label{snr4}
\end{eqnarray}
{\it provided}
\begin{eqnarray}
(r'-r)p=2r(3l-1)\pm 2l. 
\label{snr44}
\end{eqnarray}
{\it If $r, r'$ solve} 
\begin{eqnarray}
(r'-r)p=2r(3l-1)\mp 2\left(l+\frac{p-1}{3}\right) 
\label{snr41}
\end{eqnarray}
{\it then the common non-ranks are}
\begin{eqnarray}
p[p'(2n+1)+2r']\pm 4N(\frac{p}{6})=p'[p(2n+1)
+2r]\mp 4N(\frac{p'}{6}).
\label{snr42}
\end{eqnarray}
Note that, again, there are $4=2^{\nu(pp')}$ 
arithmetic progressions of common or double 
non-ranks to the primes $p', p$ in all cases. 

{\it Proof.} By substituting $p',N(p'/6)$ 
in terms of $p,N(p/6)$ and $l,$ respectively, 
it is readily verified that Eqs.~(\ref{snr1}), 
(\ref{snr11}) are equivalent, as are 
(\ref{snr12}), (\ref{snr13}), and (\ref{snr2}), 
(\ref{snr22}), and (\ref{snr23}), (\ref{snr24}), 
and (\ref{snr3}), (\ref{snr33}), and  
(\ref{snr31}), (\ref{snr32}), and (\ref{snr4}), 
(\ref{snr44}), and (\ref{snr41}), (\ref{snr42}).  
As in (i) there is a unique solution 
$(r, r')$ in all other cases as well.~$\diamond$ 

{\bf Example~3.10.} For $p=5,~p'=11$ we have $l=1$ 
and Eq.~(\ref{snr11}) becomes $5(r'-r)=2(3r+1),$ 
i.e. $r=3,~r'-r=4$ or $r'=7.$ So Eq.~(\ref{snr1})  
gives the common non-ranks 
\begin{eqnarray}
5[11\left(2(n+1)+1\right)-8]+4=11[5\left(2(n+1)
+1\right)-4]+4\cdot 2. 
\end{eqnarray} 
The other sign in Eq.~(\ref{snr11}) is $5(r'-r)
=2(3r-1)$ solved by $r=2,~r'-r=2,~r'=4$ and the 
common non-ranks are 
\begin{eqnarray}
5[11(2n+1)+8]-4=11[5(2n+1)+4]-4\cdot 2. 
\end{eqnarray} 

{\bf Theorem~3.11. (Triple non-ranks)} {\it Let 
$5\leq p<p'<p''$ (or $5\leq p<p''<p'$, or 
$5\leq p''<p<p'$) be different odd primes. Then 
each case in Theor.~3.9 of four double non-ranks 
leads to $8=2^{\nu(pp'p'')}$ triple non-ranks of 
$p, p', p''.$ At two non-ranks per prime, there 
are at most $2^3$ triple non-ranks.}  

{\it Proof.} It is based on Theor.~3.9 and 
similar for all its cases. Let's take (i) 
and substitute $2n+1\to p''(2n+1)+2\nu,~-p''<2\nu<p''$ 
in Eq.~(\ref{snr1}) which, upon dropping the term 
$p''p'p(2n+1),$ yields on the lhs 
\begin{eqnarray}
2pp'\nu+2pr'-4N(\frac{p}{6})=2p''\mu\pm 4N(\frac{p''}{6}).
\label{31}
\end{eqnarray} 
Since $(pp',p'')=1$ there is a unique residue $\nu$ 
modulo $p''$ so that the lhs of Eq.~(\ref{31}) is 
$\equiv \pm 4N(\frac{p''}{6})\pmod{p''},$ and this 
determines $\mu.$ As each sign case leads to such 
a triple non-rank solution, it is clear that there 
are $2^3$ non-ranks to $p, p', p''.~\diamond$ 
 
{\bf Example~3.12.} For $5,7,11$ triple non-rank 
progressions are obtained as follows. Starting 
from the double non-rank equations (cf. proof of 
Lemma~3.4) 
\begin{eqnarray}
5\cdot 7(2n+1)\pm 4,~5[7(2n+1)+4]+4=7[5(2n+1)+4]-4
\end{eqnarray}
replace $2n+1\to 11(2n+1)+2\nu,$ drop $5\cdot 7\cdot 
11(2n+1)$ and set the rhs to $22\mu+4N(11/6):$ 
\begin{eqnarray}
5\cdot 7\cdot 2\nu+4=11\cdot 2\mu+8.
\end{eqnarray}
Since $5\cdot 7\cdot 2-4=11\cdot 6,$ the solution 
is $\nu=1,~\mu=3.$ Putting back $5\cdot 7\cdot 
11(2n+1)$ we obtain the triple non-rank system 
\begin{eqnarray}\nonumber 
&&5\cdot 7[11(2n+1)+2]+4=7[5\cdot 11(2n+1)+10]+4\\
&&=11[5\cdot 7(2n+1)+6]+2\cdot 4. 
\end{eqnarray}
Setting the rhs to $22\mu-8$ yields the second such  
solution  
\begin{eqnarray}
35\cdot 2\nu+4=22\mu-8,~\nu=4,~\mu=13,
\end{eqnarray}
with the complete triple non-rank system 
\begin{eqnarray}\nonumber
&&5\cdot 7[11(2(n+1)+1)-6]+4=7\cdot 5[11(2(n+1)+1)-6]
+4\\&&=11[5[7\left(2(n+1)+1\right)-4]+2]-4\cdot 2.
\end{eqnarray}

{\bf Theorem~3.13. (Multiple non-ranks)} {\it Let 
$5\leq p_1<\cdots <p_m$ be $m$ different primes. 
Then there are $2^m$ arithmetic progressions of 
$m-$fold non-ranks to the primes} $p_1,\ldots,p_m.$ 

{\it Proof.} This is proved by induction on $m.$   
Theors.~3.9 and 3.11 are the $m=2, 3$ cases. If 
Theor.~3.13 is true for $m$ then for any case 
$5\leq p_{m+1}<p_1<\cdots<p_m,$ or $\ldots,$ 
$5\leq p_1<\cdots<p_{m+1},$ we substitute in 
an $m-$fold non-rank equation $2n+1\to p_{m+1}
(2n+1)+2\nu$ as in the proof of Theor.~3.11, 
again dropping the $(2n+1)\prod_1^{m+1} p_i$ 
term. Then we get   
\begin{eqnarray}\nonumber
&&p_1(p_2(\cdots(2p_m\nu+2r_m)+\cdots+2r_2)+
4N(\frac{p_1}{6})\\&&=2p_{m+1}\mu\pm 
4N(\frac{p_{m+1}}{6})
\label{m}
\end{eqnarray}  
with a unique residue $2\nu\pmod{p_{m+1}}$ so 
that the lhs of Eq.~(\ref{m}) becomes 
$\equiv 4N(\frac{p_{m+1}}{6})\pmod{p_{m+1}},$ 
which then determines $2\mu.$ In case the lhs of 
Eq.~(\ref{m}) has $p_1(\ldots)-4N(p_1/6)$ the 
argument is the same. This yields an $(m+1)-$fold 
non-rank progression since each sign in 
Eq.~(\ref{m}) gives a solution. Hence there 
are $2^{m+1}$ such non-ranks. At two non-ranks 
per prime there are at most $2^m$ non-rank 
progressions.~$\diamond$  

\section{Counting Non-Ranks}

If we subtract for case (i) in Theor.~3.9, say, 
the four common non-rank progressions this leaves in  
${\cal A}_{p'}^-=\{p'(2n+1)\pm 4\frac{p'+1}{6}: 
n\geq 0\}$ the following progressions $p'p(2n+1)\pm 
4\frac{p'+1}{6},\ldots, p'[(2n+1)p+2r_1]+4\frac{p'
+1}{6},\ldots, p'[(2n+1)p+2r_2]-4\frac{p'+1}{6},
\ldots, p'[(2n+1)p+2r_3]+4\frac{p'+1}{6},\ldots, 
p'[(2n+1)p+2r_4]-4\frac{p'+1}{6},\ldots, 
p'(2n+1)p\pm 4\frac{p'+1}{6}.$ 
 
We summarize this as follows. 

{\bf Lemma~4.1.} {\it $p'>p\geq 5$ be prime.
Removing the common non-ranks of $p',p$ 
from the set of all non-ranks of $p'$ 
leaves arithmetic progressions of the 
form $p'(2n+1)p+2l;~n\geq 0$ where $l>0$ are 
given non-negative integers.}  

{\bf Proposition~4.2.} {\it Let $p>5$ be prime. 
Then the set of non-ranks to parent prime $p,~
{\cal A}_p,$ is made up of arithmetic 
progressions $L(p)(2n+1)+2a,~n\geq 0$ with $L(p)=
\prod_{5\leq p'\leq p}p,~p'$ prime and $a>0$ 
given integers.}  

{\it Proof.} Let $p=6m\pm 1.$ We start from the 
set ${\cal A}_p^\pm=\{p(2n+1)\pm 4N(\frac{p}{6})
>0: n=0,1,2,\ldots\}.$ Removing the non-ranks 
common to $p$ and $5$ by Lemma~4.1 leaves 
arithmetic progressions of the form $5p(2n+1)
+2l,~n\geq 0$ where $l>0$ are given integers. 
Continuing this process to the largest 
prime $p'<p$ leaves in ${\cal A}_p$ arithmetic 
progressions of the form $L(p)(2n+1)+2a,~n\geq 0$ 
with $L(p)=\prod_{5\leq p'\leq p}p'$ and $a>0$ 
a sequence of given integers independent of 
$n.~\diamond$         

{\bf Proposition~4.3.} {\it Let $p\geq 
p'\geq 5$ be primes and $G(p)$ the number 
of non-ranks $L(p)(2n+1)+2a\in {\cal A}_p$ 
over one period $L(p)$ corresponding to 
arithmetic progressions $L(p)(2n+1)+2a\in 
{\cal A}_p.$ Then} $G(p)=\prod_{5\leq p'<p}
(p'-2).$  

Note that $G(p)<L(p)$ both increase 
monotonically as $p\to\infty.$ 

{\it Proof.} In order to determine $G(p)$ 
we have to eliminate all non-ranks of 
primes $5\leq p'<p$ from ${\cal A}_p.$ 
As in Lemma~3.4 we start by subtracting 
the fraction $2/5$ from the interval 
$1\leq a\leq L(p)$ of length $L(p),$ then 
$2/7$ for $p'=7$ and so on for all $p'<p.$ 
The factor of $2$ is due to the symmetry 
of non-ranks around each multiple of $p'$ 
according to Lemma~2.5. This leaves $p
\prod_{5\leq p'<p}(p'-2)/2$ odd numbers.  
The fraction $2/p$ of these are the 
non-ranks to parent prime $p.~\diamond$   

Prop.~4.3 implies that the fraction of 
non-ranks related to a prime $p$ in the 
interval occupied by ${\cal A}_p,$  
\begin{eqnarray}
q(p)=\frac{G(p)}{L(p)}=\frac{1}{p}
\prod_{5\leq p'<p}\frac{p'-2}{p'},
\label{qp}
\end{eqnarray}
where $p'$ is prime, decreases 
monotonically as $p$ goes up.  

{\bf Definition~4.4.} Let $p\geq p'\geq 5$ be 
prime. The supergroup ${\cal S}_p=\bigcup_{
p'\leq p}\\{\cal A}_{p'}$ contains the sets 
of non-ranks corresponding to arithmetic 
non-rank progressions $2a+L(p')(2n+1)$ of 
all ${\cal A}_{p'},~p'\leq p.$   

Thus, each supergroup ${\cal S}_p$ contains 
nested sets of non-ranks related to primes 
$5\leq p'\leq p.$

Let us now count prime numbers from $p_1=2$ 
on. 

{\bf Proposition~4.5.} {\it Let $p_j\geq 5$ 
be the $j$th prime. (i) Then the number of 
non-ranks $a\in {\cal A}_{p_i}$ 
corresponding to arithmetic progressions 
related to a prime $5\leq p_i<p_j,$}
\begin{eqnarray}
G(p_i)=\frac{L(p_j)}{L(p_i)}G(p_j)=\frac{L(p_j)}
{p_i}\prod_{5\leq p<p_i}\frac{p-2}{p}=q(p_i)L(p_j),
\end{eqnarray}
{\it where $p$ is prime, monotonically decreases 
as $p_i$ goes up. (ii) The number of non-ranks in 
a supergroup ${\cal S}_{p_j}$ over one period 
$L(p_j)$ is}
\begin{eqnarray}
S(p_j)=L(p_j)\sum_{5\leq p\leq p_j}q(p)=\frac{1}{2}
L(p_j)\left(1-\prod_{5\leq p\leq p_j}\frac{p-2}{p}
\right).\label{lsp} 
\end{eqnarray}
{\it (iii) The fraction of non-ranks of their 
arithmetic progressions in the (first) interval 
$[1,L(p_j)]$ occupied by the supergroup} 
${\cal S}_{p_j},$
\begin{eqnarray}
Q(p_j)=\frac{S(p_j)}{L(p_j)}=\sum_{5\leq p\leq p_j}
q(p)=\frac{1}{2}[1-\prod_{5\leq p\leq p_j}\frac{p-2}{p}],
\label{qp1}
\end{eqnarray}
{\it increases monotonically as $p_j$ goes up.} 

{\it Proof.} (i) follows from Prop.~4.3 
and Eq.~(\ref{qp}). (ii) and (iii) are 
equivalent and are proved by induction as follows, 
using Def.~4.4 in conjunction with Eq.~(\ref{qp}).  

From Eq.~(\ref{qp}) we get $q_3=2/p_3$ which is 
the case $j=3,~p_j=5$ of Eq.~(\ref{qp1}). Assuming 
Eq.~(\ref{qp1}) for $p_j,$ we add $q_{j+1}$ of 
Eq.~(\ref{qp}) and obtain
\begin{eqnarray}\nonumber
\sum_{i=3}^{j+1}q(p_i)&=&\frac{1}{2}-\frac{1}{2}
\prod_{i=3}^j\frac{p_i-2}{p_i}+\frac{1}{p_{j+1}}
\prod_{i=3}^j\frac{p_i-2}{p_i}\\&=&\frac{1}{2}-
\frac{1}{2}\prod_{i=3}^{j+1}\frac{p_i-2}{p_i}.
\end{eqnarray}
The extra factor $0<(p_{j+1}-2)/p_{j+1}<1$ shows 
that $q(p_j), x(p_j)$ in Eq.~(\ref{xp}) decrease 
monotonically as $p_j\to p_{j+1}$ while $Q(p_j)$ 
increases as $j\to\infty.~\diamond$

{\bf Definition~4.6.} Since $L(p)>S(p),$ there 
is a set ${\cal R}_p$ of {\it remnants} $r\in 
[1, L(p)]$ such that $r\not\in{\cal S}_p.$   

{\bf Lemma~4.7.} {\it (i) The number $R(p_j)$ of 
remnants in a supergroup, ${\cal S}_{p_j},$ is}   
\begin{eqnarray}
R(p_j)=\frac{1}{2}L(p_j)-S(p_j)=L(p_j)(\frac{1}{2}
-Q(p_j))=\frac{1}{2}\prod_{5\leq p\leq p_j}(p-2)
=\frac{1}{2}G(p_{j+1}).
\label{rps}
\end{eqnarray} 

{\it (ii) The fraction of remnants in} 
${\cal S}_{p_j},$   
\begin{eqnarray}
x(p_j)&=&\frac{R(p_j)}{L(p_j)}=\frac{1}{2}-Q(p_j)
=\frac{1}{2}\prod_{5\leq p\leq p_j}\frac{p-2}{p},
\label{xp}
\end{eqnarray}
{\it where $p$ is prime, decreases monotonically 
as} $p_j\to\infty.$ 

{\it Proof.} (i) follows from Def.~4.6 in 
conjunction with Eq.~(\ref{lsp}) and (ii) 
from Eq.~(\ref{rps}). Eq.~(\ref{rps}) follows 
from Eq.~(\ref{qp1}).~$\diamond$
       
\section{Remnants and Twin Ranks}

When all primes $5\leq p\leq p_j$ and appropriate 
nonnegative integers $n$ are used in Lemma~2.5 
one will find all non-ranks $2k+1<M(j+1)\equiv 
(p_{j+1}^2-2^2)/3.$ By subtracting these non-ranks 
from the set of odd positive integers $N\leq M(j+1)$ 
all and only twin ranks $t<M(j+1)$ are left among 
the remnants. If a non-rank $2k+1$ is left then 
$3(2k+1)\pm 2$ must have prime divisors that are 
$>p_j$ according to Lemma~2.5, which is 
impossible. All $t<M(j+1)=(p_{j+1}^2-4)/3$ in a 
remnant ${\cal R}_{p_j}$ of a supergroup 
${\cal S}_{p_j}$ are twin-4 ranks. 

{\bf Proposition~5.1.} {\it Let $p_j$ be the $j$th  
prime number and $L(p_j)(2n+1)+a_i^{(j)}$ be the 
arithmetic progressions that contain the non-ranks 
$a_i^{(j)}\in{\cal A}_{p_j}$ to parent prime $p_j.$ 
Let $3[L(p_j)(2n+1)+c_i^{(j)}]\pm 2$ be the 
arithmetic progressions that contain the twin-4 
primes with $c_i^{(j)}\in{\cal C}_{p_j}.$   

(i) The set of constants $c_i^{(j+1)}$ of 
arithmetic progressions containing the twin-4 
ranks from the next supergroup 
${\cal S}_{p_{j+1}}$ is} 
\begin{eqnarray}\nonumber
&&{\cal C}_{p_{j+1}}=\{3[L(p_j)(p_{j+1}(2n+1)+l)
+c_i^{(j)}]\pm 2\}\\&&-\{3[L(p_j)(p_{j+1}(2n+1)+l')
+a_{i'}^{(j)}]\pm 2\}. 
\label{rec}
\end{eqnarray}
{\it If there are positive integers $0\leq l,l'<p_{j+1},$ 
a non-rank $a_{i'}^{(j)}\in {\cal A}_{p_j}$ and a 
constant $c_i^{(j)}\in {\cal C}_{p_j}$ satisfying} 
\begin{eqnarray}
L(p_j)l+c_i^{(j)}=L(p_j)l'+a_{i'}^{(j)}, 
\end{eqnarray}
{\it then} 
\begin{eqnarray}
L(p_j)l+c_i^{(j)}\not\in {\cal C}_{p_{j+1}},
\end{eqnarray}
{\it else}
\begin{eqnarray} 
c_{i,l}^{(j+1)}=L(p_j)l+c_i^{(j)}\in {\cal C}_{p_{j+1}}. 
\end{eqnarray}   
Prop.~5.1 is the inductive step completing the 
practical sieve construction for ordinary twin 
primes. Props.~3.3, 3.5 and Lemma~3.4 are the 
initial steps.  

{\it Proof.} Replacing in (i) $2n+1\to p_{j+1}(2n+1)
+l,~l=0,1,2,\ldots, p_{j+1}-1$ and subtracting the 
resulting sets from each other, we obtain 
(i).$~\diamond$  

For $p_3=5,$ Prop.~5.1 is Prop.~3.3, for $p_4=7$ 
it is Prop.~3.5. Clearly, at the start of the 
$c$ for $p_4=7$ the previous values for $p_3=5$ 
are repeated, but this pattern does not continue. 
 
Twin-4 ranks are located among the remnants 
${\cal R}_p$ for any prime $p\geq 5$. Our 
goal is to develop a Legendre-type sum for the 
number $R$ of twin-4 ranks.   

{\bf Theorem~5.2.} {\it Let $R_0$ be the number 
of remnants of the supergroup ${\cal S}_{p_j},$ 
where $p_j$ is the $j$th prime number and 
$M(j+1)=[p_{j+1}^2-4]/3$. Then the number 
$R=\pi_2(3L(p_j)+2)/2$ of twin-4 ranks 
within the remnants of the supergroup 
${\cal S}_{p_j}$ is given by} 
\begin{eqnarray}
R=R_0+\sum_{p_j<n}\mu(n)2^{\nu(n)}\bigg[
\frac{L(p_j)-M(j+1)-1}{2n}\bigg].
\label{tr} 
\end{eqnarray}

Here $L(p_j)=\prod_{5\leq p\leq p_j}p,$  
$R_0=\frac{1}{2}\prod_{5\leq p\leq p_j}(p-2)$ 
with $p$ prime, and $n$ runs through all products 
of primes $p_j<p\leq (6L(p_j)+1)/4$. The upper 
limit $(6L(p_j)+1)/4$ comes about because $4N(p/6)$ 
is the lowest non-rank of the prime number $p$ 
according to Lemma~2.2.  

The argument of the twin-prime counting function 
$\pi_2$ is $3L(p_j)+2$ because, if $L(p_j)$ is the 
last twin-4 rank of the interval $[1, L(p_j)],$ then 
$3L(p_j)\pm 2$ are the corresponding twin-4 primes.  

{\it Proof.} According to Prop.~4.5 the 
supergroup ${\cal S}_{p_j}$ has $S(p_j)=\frac{1}
{2}L(p_j)\cdot\left(1-\prod_{5\leq p\leq 
p_j}\frac{p-2}{p}\right)$ non-ranks. Subtracting 
these from the interval $[1, L(p_j)]$ that the 
supergroup occupies gives $R_0=\frac{1}{2}\prod_{
5\leq p\leq p_j}(p-2)$ for the number of remnants 
which include twin-4 ranks and non-ranks to primes 
$p_j<p\leq (6L(p_j)+1)/4.$ The latter are 
\begin{eqnarray}
M(j+1)<p(2n+1)\pm 4N(\frac{p}{6})\leq L(p_j),~
M(j+1)=\frac{p_{j+1}^2-4}{3}, 
\end{eqnarray} 
or 
\begin{eqnarray}
0<n\leq \frac{L(p_j)-M(j+1)-1}{2p}, 
\end{eqnarray}
which have to be subtracted from the remnants 
to leave just twin-4 ranks. Correcting for 
double counting of common non-ranks to two 
primes using Theor.~3.9, of triple non-ranks 
using Theor.~3.11 and multiple non-ranks 
using Theor.~3.13 we obtain 
\begin{eqnarray}\nonumber
R&=&R_0-2\sum_{p_j<p\leq (6L(p_j)+1)/4}\bigg[\frac{L(p_j)
-M(j+1)-1}{2p}\bigg]\\&+&4\sum_{p_j<p<p'\leq (6L(p_j)+1)/4}
\bigg[\frac{L(p_j)-M(j+1)-1}{2pp'}\bigg]\mp\cdots, 
\label{tr1}
\end{eqnarray} 
where $[x]$ is the integer part of $x$ as 
usual. Note that the arithmetic details (functions 
of $r,r'$ in double non-ranks in Theor.~3.9, etc., 
that do not depend on $n$) do not affect their 
counting in Eq.~(\ref{tr1}) because they always 
add to $p'p(2n+1), pp'p''(2n+1), \ldots .$ 
Equation~(\ref{tr1}) is equivalent to 
Eq.~(\ref{tr}).~$\diamond$ 

{\bf Definition~5.3.} Decomposing the floor function 
$[x]=x-\{x\}$ in Eq.~(\ref{tr}) allows writing 
$R=R_M+R_E$ in terms of a main and error term
\begin{eqnarray}\nonumber
R_M&=&R_0+\sum_{p_j<n}\mu(n)2^{\nu(n)}\frac{L(p_j)
-M(j+1)-1}{2n},\\
R_E&=&-\sum_{p_j<n}\mu(n)2^{\nu(n)}\{\frac{L(p_j)
-M(j+1)-1}{2n}\}.
\label{tr2}
\end{eqnarray}

{\bf Theorem~5.4.} {\it The main term $R_M$ in 
Eq.~(\ref{tr}) satisfies}  
\begin{eqnarray}\nonumber
&&R_M=\frac{1}{2}L(p_j)\prod_{5\leq p\leq 
(6L(p_j)+1)/4}\left(1-\frac{2}{p}\right)\\&&
+\frac{1}{2}M(j+1)[1-\prod_{p_j<p\leq 
(6L(p_j)+1)/4}(1-\frac{2}{p})]. 
\label{tr3}
\end{eqnarray}

{\it Proof.} Expanding the product  
\begin{eqnarray}
L(p_j)\prod_{5\leq p\leq p_j}(1-\frac{2}{p})
\end{eqnarray}
and combining corresponding sums in Eq.~(\ref{tr2})
\begin{eqnarray}
-\sum_{5\leq p\leq p_j}\frac{1}{p}-\sum_{p_j<p\leq 
(6L(p_j)+1)/4}\frac{1}{p}=-\sum_{5\leq p\leq 
(6L(p_j)+1)/4}\frac{1}{p},\ldots 
\end{eqnarray} 
just shifts the upper limit of the primes in 
the product $\prod_p(1-2/p)$ from $p_j$ to 
$(6L(p_j)+1)/4,$ so that we obtain 
Eq.~(\ref{tr3}). This involves considerable 
cancellations collapsing $R_0$ to the correct 
magnitude of $R_M$.~$\diamond$

{\bf Theorem~5.5.} {\it The main term $R_M$ 
obeys the asymptotic law}   
\begin{eqnarray}
R_M\sim \frac{2c_2 e^{-2\gamma}3L(p_j)}
{\log^2((3L(p_j)+0.5)/2)},~p_j\to\infty. 
\end{eqnarray}

{\it Proof.} This follows as Theor.~5.8 in 
Ref.~\cite{adhjw}.~$\diamond$  

\section{Summary and Discussion} 

The cousin prime sieve is specifically designed 
for prime twins at distance $4,$ often called 
cousin primes.   

Accurate counting of non-rank sets require 
the infinite, but sparse set of odd 'primorials' 
$\{3L(p_j)=\prod_{2<p\leq p_j}p\}$ much like in 
other sieves when applied to primes. The twin-4 
primes are not directly sieved, rather twin-4 
ranks $2m+1$ are with $3(2m+1)\pm 2$ both prime. 
All other odd numbers ($\geq 9$) are non-ranks. 
Primes serve to organize and classify non-ranks 
in arithmetic progressions with equal distances 
(periods) that are primes or products of them 
leading to the (odd) primorials.       

The coefficient $4c_2 e^{-2\gamma}
\approx 0.8324267$ of the asymptotic form of 
the main term is the same as for ordinary twins, 
despite the differences in the arithmetic of the 
pair sieves. Just as for ordinary twin primes, 
the resolution of the parity problem allows 
replacing the need for a lower bound on $R,$ 
or $\pi_2$ at primorials, by an upper bound 
for the error term $R_E.$


\end{document}